\documentclass[a4paper,11pt]{article}
\title{Equilateral Sets in infinite dimensional\\ Banach Spaces}
\author{S.K.Mercourakis and G.Vassiliadis}
\date{}
\usepackage{amsmath}
\usepackage{amsfonts}
\usepackage{amssymb}
\usepackage{amscd}
\usepackage{amsthm}
\usepackage{makeidx}
\usepackage{euscript}
\theoremstyle{plain}
\newtheorem{theo}{Theorem}

\newtheorem{prop}{Proposition}
\newtheorem{cor}{Corollary}
\theoremstyle{definition}
\newtheorem{defn}{Definition}
\newtheorem*{claim}{\underline{Claim}}
\newcommand{\beg}{\begin{proof}[Proof]}
\newcommand{\eq}{equilateral $\,$}
\newcommand{\ind}{infinite dimensional Banach $\,$}
\newcommand{\ban}{Banach space $\,$}

\begin{document}
\maketitle

\begin{abstract}
We show that every \ban X containing an isomorphic copy of $c_0$
has an infinite \eq set and also that if X has a bounded
biorthogonal system of size $\alpha$ then it can be renormed so as
to admit an \eq set of equal size.
\end{abstract}


\section*{Introduction}
A subset S of a metric space $(M,d)$ is said to be \eq if there is
a $\lambda>0$ such that for $x,y \in S,x \neq y$ we have
$d(x,y)=\lambda$; we also call S a $\lambda$-\eq set. Equilateral
sets have been studied mainly in finite dimensional spaces (see
\cite{Pe},\cite{S1},\cite{S2},\cite{SV}).

Our aim in this note is the study of \eq sets in infinite
dimensions. We first prove (improving results of K.J. Swanepoel)
that any \ind space has an equivalent norm arbitrarily close to
the original one admitting an infinite \eq set (Th.1). Then we
prove that every \ban containing an isomorphic copy of $c_0$
admits an infinite \eq set (Th.2). We also introduce a notion of
antipodal sets (Defs. 1 and 2) which yields that a \ban containing
a bounded biorthogonal system of size $\alpha$ can be renormed so
that in the new (equivalent) norm it has an \eq set of equal size
(Th.3). These results generalize results of Petty \cite{Pe} and
Swanepoel \cite{S1},\cite{S2}.

If X is any (real) \ban then $B_X$ denotes its closed unit ball.
The Banach-Mazur distance between two isomorphic Banach spaces X
and Y is\\ $d(X,Y)=\inf\{\|T\| \cdot \|T^{-1}\|$  where $T:X
\rightarrow Y $ is a linear isomorphism $\}$.

A sequence $(e_n)$ in a Banach space X is said to be spreading, if
for any sequence $0<p_1<p_2< \dots<p_N$ of integers and any
sequence $\alpha_1,\alpha_2,\dots,\alpha_N$ of scalars we have
$\|\sum_{k=1}^N \alpha_k e_k\|=\|\sum_{k=1}^N\alpha_k e_{p_k}\|$.
It is clear that a non-constant spreading sequence is equilateral.
About the theory of spreading models we refer the reader to
\cite{BL} and \cite{Alkal}.

\underline{Acknowledgements.} We wish to thank the referee for his
valuable suggestions and remarks. In particular, Theorem 1 in this
general form is due to the referee. Initially we had proved a
weaker result.

\section*{Equilateral and antipodal sets in infinite dimensions}

The question whether any \ind space contains an infinite \eq set
has been answered in the negative by P.Terenzi in \cite{Ter}. The
question we are concerned with in this note is under which
conditions an \ind space contains an infinite \eq set.

K.J.Swanepoel has proved that every \ind space may be renormed so
as to contain an infinite \eq set and also that if the space is
uniformly convexifiable then we may choose the new norm to be
arbitrarily close to the original norm (see \cite{S1},\cite{S2}).

P.Brass \cite{Bra} and B.V. Dekster \cite{Dek} have proved (both
using Dvoretzky's Theorem) that each infinite dimensional Banach
space contains arbitrarily large (finite) equilateral sets (see
also \cite{S2}).

Our first result improves the results of Swanepoel mentioned
above.

\begin{theo}
Let $(X,\|\cdot\|)$ be an infinite dimensional Banach space. Then
for every $\varepsilon >0$ there exists an equivalent norm $|||
\cdot |||$ on X such that

\begin{enumerate}
\item $d((X,\| \cdot \|),(X,||| \cdot |||))\le 1+\varepsilon$
\item $(X,||| \cdot |||)$ admits an infinite \eq set.
\end{enumerate}
\end{theo}

\beg It is enough to find a closed subspace Z of X and an
equivalent norm $|||\cdot|||$ on Z, so that conditions 1. and 2.
are satisfied for Z (see Remark 1(1) below).

Assume first that X contains an isomorphic copy, say Y of
$\ell_1$. By the classical non distortion property of $\ell_1$
(\cite{LT}, Prop. 2.e.3) for every $\varepsilon>0$ there is a
normalized sequence $(y_n)$ in $(Y,\| \cdot \|)$ such that for
every sequence $(\alpha_n) \in c_{00}$ (the space of eventually
zero real sequences) we have
\[ \frac{1}{1+\varepsilon} \sum_{n=1}^\infty |\alpha_n| \le \|
\sum_{n=1}^\infty \alpha_n y_n \| \le \sum_{n=1}^\infty
|\alpha_n|.\]

Let $Z=\overline{<y_n>}$ (i.e. the closed linear span of $(y_n)$)
and for $y=\sum_{n=1}^\infty \alpha_n y_n \in Z$, set
$|||y|||=\sum_{n=1}^\infty |\alpha_n|$. Clearly the space $(Z,|||
\cdot |||)$  is isometric to $\ell_1$ and $d((Z,\| \cdot
\|),(Z,||| \cdot |||))\le 1+\varepsilon$.

Now assume that $\ell_1 \not \hookrightarrow X$, so there exists a
normalized weakly null (basic) sequence $(x_n)$ in X with
spreading model $(e_n)$, which is a normalized unconditional basic
sequence with suppression constant $K_s=1$ (see \cite{Alkal} pp.
275-279). Given $\varepsilon>0$ and $m \in \Bbb{N}$ ($m \ge 2$) it
is enough to produce a subsequence of $(x_n)$ (still denoted
$(x_n)$) and a norm $||| \cdot|||$ on
the span $<x_n>$ of $(x_n)$ satisfying the following:\\
(i) $\frac{1}{1+\varepsilon} \|x\| \le |||x||| \le
\frac{1}{(1-\varepsilon)^2} \|x\|$
for $x=\sum \alpha_n x_n, (\alpha_n) \in c_{00}$ and\\
(ii)$|||\sum_{n \in F} \alpha_n x_n |||=\|\sum_{n \in F} \alpha_n
e_n \|$ for all $F \subseteq \Bbb{N}$ with $|F|=m$.

Indeed, by passing to a subsequence we may assume that\\
(a) $(1-\varepsilon) \|\sum_{n \in F} \alpha_n e_n\| \le \|
\sum_{n \in F} \alpha_n x_n\| \le (1+\varepsilon) \|\sum_{n \in F}
\alpha_n e_n\|$ for all $(\alpha_n) \in c_{00}$ and $F \subseteq
\Bbb{N}$ with $|F|=m$ and\\
(b) $\|\sum_{n=1}^\infty \alpha_n x_n\| \ge (1-\varepsilon)
\|\sum_{n \in F} \alpha_n x_n\|$ for all $(\alpha_n) \in c_{00}$
and $|F|=m$, by using Schreier unconditionality and keeping in
mind that $m$ is fixed (see \cite{Od}).

Let $x=\sum \alpha_n x_n,(\alpha_n) \in c_{00}$; we set
\[ |||x|||=max\{\frac{1}{1+\varepsilon} ||x||,\, \sup_{|F|=m} \|\sum_{n \in F} \alpha_n
e_n\| \}. \]

It is obvious that $\frac{1}{1+\varepsilon} \|x\| \le |||x|||$ and
since
\[ \|\sum_{n \in F} \alpha_n e_n\| \le \frac{1}{1-\varepsilon} \|\sum_{n \in F} \alpha_n
x_n\| \le \frac{1}{(1-\varepsilon)^2} \|x\|\] by (a) and (b) we
get (i).

Let now $y=\sum_{n \in F_0} \alpha_n x_n, |F_0|=m$, then (a)
implies that \[\|\sum_{n \in F_0} \alpha_n e_n\| \ge
\frac{1}{1+\varepsilon} \|y\|\]  so since the unconditional
sequence $(e_n)$ has suppression constant $K_s=1$ it follows
$|||y|||=\|\sum_{n \in F_0} \alpha_n e_n\|$ and hence we get (ii).

Set $Z=\overline{<x_n>}$; since the sequence $(e_n)$ is spreading,
by (ii) we get the conclusion.
\end{proof}

\noindent \textbf{Remarks 1} (1) Let $(X,\|\cdot\|)$ be a normed
space, Z a linear subspace of X and $|||\cdot|||$ an equivalent
norm on Z such that \[c_1 |||x||| \le \|x\| \le c_2 |||x||| \,\,
\forall x \in Z \,\,\, (1).\] Then the norm $|||\cdot|||$ of Z can
be extended on X (using the Hahn-Banach Theorem) so that (1) is
satisfied. Indeed, for $x \in X$ set
$|||x|||_1=\sup|\tilde{f}(x)|$, where the supremum is taken over
all $f \in Z^*$ with $|||f||| \le 1$ and $\tilde{f}$ is a
Hahn-Banach extension of $f$ on X with $||f||=||\tilde{f}||$. Then
take $|||x|||=\max\{|||x|||_1,\frac{1}{c_2} \|x\|\}$ and it is
easy to see that (1) is satisfied for all $x \in X$.

(2) By a result of Rosenthal (\cite{Ros}) the sequence $(x_n)
\subseteq X$ in the proof of Th.1 can be chosen so that the
spreading sequence $(e_n)$ is 1-unconditional (i.e. with
unconditional constant $K_u=1$). Therefore the proof gives more
than equilateral; in particular we get that $|||x_n \pm
x_m|||=|||x_1-x_2|||>0$ for all $n,m \in \Bbb{N}$ with $n \neq m$.
More generally $|||\sum_{i=1}^m a_i x_{n_i}|||=|||\sum_{i=1}^m b_i x_{k_i}|||$
whenever $n_1<n_2<\cdots <n_m,\, k_1<k_2<\cdots <k_m$ and $|a_i|=|b_i|$ 
for $i \le m$.

\begin{theo}
Every \ban X containing an isomorphic copy of $c_0$ admits an
infinite equilateral set.
\end{theo}

\beg We shall use the non-distortion property of $c_0$ and the
following generalization of Theorem B of \cite{SV}, with similar
proof (see also \cite{Bra}).

\begin{claim}
Let $\|\cdot\|$ be an equivalent norm on $c_0$ with Banach-Mazur
distance at most $\frac{3}{2}$ from the original norm
$\|\cdot\|_\infty$ of $c_0$. Then $(c_0,\|\cdot\|)$ admits an
infinite \eq set.
\end{claim}

\underline{Proof of the Claim:} We may assume that $\|x\| \le
\|x\|_\infty \le \frac{3}{2} \|x\|$ for $x \in c_0$. Let
$I=\{(n,m):n,m \in \Bbb{N}$ and $n<m\}$; denote by K the compact
cube $[0,\frac{1}{2}]^I$.

For $\varepsilon=(\varepsilon_{(n,m)})
\in K$ we set: $p_1(\varepsilon)=(-1,0,\dots)$ and\\
$p_n(\varepsilon)=(\varepsilon_{(1,n)},\varepsilon_{(2,n)},\dots,
\varepsilon_{(n-1,n)},-1,0,\dots)$ for $n \ge 2$. Observe that for
$n<m$ we have
\[\|p_n(\varepsilon)-p_m(\varepsilon)\|_\infty=1+\varepsilon_{(n,m)}.\]
We define a function $\varphi:K \rightarrow K$ by the rule
$\varphi_{(n,m)}(\varepsilon)=1+\varepsilon_{(n,m)}-
\|p_n(\varepsilon)-p_m(\varepsilon)\|, \, (n,m) \in I, \,
\varepsilon \in K$. Note that $\varphi_{(n,m)}(\varepsilon) \ge
1+\varepsilon_{(n,m)}-
\|p_n(\varepsilon)-p_m(\varepsilon)\|_\infty=0$ and
$\varphi_{(n,m)}(\varepsilon) \le
1+\varepsilon_{(n,m)}-\frac{2}{3}
\|p_n(\varepsilon)-p_m(\varepsilon)\|_\infty=\frac{1}{3}
(1+\varepsilon_{(n,m)}) \le \frac{1}{2}$, so $\varphi$ is well
defined. Since each coordinate function $\varphi_{(n,m)}$ is
continuous, $\varphi$ is also continuous. Hence by a classical
result of Schauder $\varphi$ has a fixed point
$\varepsilon'=(\varepsilon'_{(n,m)}) \in K$; that is
$\varphi(\varepsilon')=\varepsilon'$, which implies that
$\|p_n(\varepsilon')-p_m(\varepsilon')\|=1$ for $n<m$. Therefore
the set $\{p_n(\varepsilon'):n \in \Bbb{N}\}$ is \eq in
$(c_0,\|\cdot\|)$ and the Claim holds.

Denote by $\|\cdot\|$ the norm on X and let Y be a subspace of X
isomorphic to $c_0$. By the non-distortion property of
$(c_0,\|\cdot\|_\infty)$ there is a subspace Z of Y (isomorphic to
$c_0$) such that $d((Z,\|\cdot\|),(c_0,\|\cdot\|_\infty)) \le
\frac{3}{2}$ (see also Th.1). It follows immediately from the
Claim that the space $(Z,\|\cdot\|)$ admits an infinite \eq set.
\end{proof}

\underline{Notes :} (1) $c_0$ cannot be replaced by $\ell_1$ in Theorem 2.
Indeed Terenzi's example gives for every $\epsilon>0$ a 
$(1+\epsilon)$-renorming of $\ell_1$ so as to not contain an 
infinite equilateral set.

(2) Let K be an infinite compact Hausdorff space. Since (as it is well known)\\
$(C(K),||\cdot||_\infty)$ contains an isometric copy of $c_0$, its unit
sphere contains a sequence $(x_n)$ with $||x_n-x_m||_\infty=2$ for $n \neq m$.
So it seems natural to ask how large a 2-equilateral subset of the
unit sphere of $C(K)$ can be, assuming further that K is (compact) 
non-metrizable. In "most" cases one can prove that such a set is uncountable,
but the general case is open for us.

In the following definition we generalize a concept coming from
finite dimensions to infinite dimensional spaces.

\begin{defn}
Let $(X,\| \cdot \|)$ be a normed space. A subset S of X is said
to be antipodal if for every $x,y \in S$ with $x \neq y$ there
exists $f \in X^*$ such that $f(x)<f(y)$ and $f(x) \le f(z) \le
f(y) \, \forall z \in S$. That is for every $x,y \in S$ with $x
\neq y$ there exist closed distinct parallel support hyperplanes
$P(=\{z \in X:f(z)=f(x)\})$ and $Q(=\{z \in X:f(z)=f(y)\})$ with
$x \in P$ and $y \in Q$.
\end{defn}

\noindent \textbf{Remarks 2} (1) If X is a finite dimensional real
vector space then the concept of antipodality coincides with the
classical one.

(2) It is well known by a result of Danzer and Gr\"{u}nbaum
\cite{DG} that the maximum cardinality of an antipodal set in
$\Bbb{R}^n$ is $2^n$ and this is attained only if the points of
the antipodal set are the vertices of an n-dimensional
parallelotope. A typical example of such a set is the unit ball B
of $\ell_\infty^n$; the vertices of B are: (the extreme points of
B) $\{(\varepsilon_1,\dots,\varepsilon_n):\varepsilon_i=\pm1,\,
i=1,2,\dots,n\}$.

Let X be a Banach space. A family $\{(x_\gamma,x_\gamma^*),\gamma
\in \Gamma \}$ of pairs in $X \times X^*$ is called a biorthogonal
system, if $x_\beta^*(x_\alpha)=\delta_{\alpha\beta}$, where
$\delta_{\alpha\beta}$ is the Kronecker $\delta$, for all
$\alpha,\beta \in \Gamma$. A family $\{x_\gamma:\gamma \in
\Gamma\}$ in X is called a minimal system, if there exists a
family $\{x_\gamma^*:\gamma \in \Gamma\}$ in $X^*$ such that
$\{(x_\gamma,x_\gamma^*),\gamma \in \Gamma \}$ is a biorthogonal
system.

\begin{prop}
Every minimal system in a \ban is antipodal.
\end{prop}

\beg Let $\{x_\gamma:\gamma \in \Gamma\}$ be a minimal system in
X, hence there exists $\{x_\gamma^*:\gamma \in \Gamma\} \subseteq
X^*$ such that the family $\{(x_\gamma,x_\gamma^*),\gamma \in
\Gamma \}$ is a biorthogonal system. Let $\gamma_1,\gamma_2 \in
\Gamma$ with $\gamma_1 \neq \gamma_2$, then we have
\[0=x_{\gamma_1}^*(x_{\gamma_2}) \le x_{\gamma_1}^*(x_{\gamma})
\le x_{\gamma_1}^*(x_{\gamma_1})=1 \,\,\, \forall \gamma \in
\Gamma.\]
It follows immediately that $\{x_\gamma:\gamma \in
\Gamma\}$ is an antipodal set in X.
\end{proof}

The following result generalizes a result of Petty with
essentially the same proof (\cite{Pe}, Th.1).

\begin{prop}
Let S be an \eq set in a normed space $(X,\| \cdot \|)$, then S is
antipodal.
\end{prop}

\beg Let $x,y \in S,x \neq y$. Suppose that S is a $\lambda$-\eq
set. By the Hahn-Banach theorem there is an $f \in X^*, \|f \|=1$
such that
\[f(y-x)=\|y-x \|=\lambda>0.\] Then $f(x)<f(y)$ and $f(y)=sup \{f(z):z \in
B(x,\lambda) \}$. So f is a support functional of the ball
$B(x,\lambda)$ through $y$ and $f(z) \leq f(y) \, \forall z \in
S$. Also if $g=-f$ we have \[g(x-y)=f(y-x)=\|y-x \|>0\] and $\|g
\|=1$, so similarly g is a support functional of the ball
$B(y,\lambda)$ through $x$ and $g(z) \leq g(x) \, \forall z \in
S$. Hence $f(x) \leq f(z) \leq f(y) \, \forall z \in S$ and the
set S is antipodal.
\end{proof}

Petty has also proved that if S is an antipodal set in a finite
dimensional real vector space X, then there exists a norm $\|
\cdot \|$ on X such that S is \eq in $(X,\| \cdot \|)$ (\cite{Pe},
Th.2). In order to generalize this result in infinite dimensions
we shall need a strengthening of the concept of antipodal set
introduced in Definition 1.

\begin{defn}
Let $(X,\| \cdot \|)$ be a normed space. We call an antipodal
subset S of X (cf. Def.1) \underline{bounded} and
\underline{separated}, if there are positive constants $c_1,c_2$
and d such that
\begin{enumerate}
\item $\|x\| \le c_1, \forall x \in S$ and
\item for every $x,y \in S$ with $x \neq y$ there is an $f \in
X^*$ with $\|f\| \le c_2$, such that $0<d \leq f(y)-f(x)$ and
$f(x) \leq f(z) \leq f(y) \, \forall z \in S$.
\end{enumerate}
\end{defn}

\noindent \textbf{Remarks 3} (1) Let S be a bounded and separated
antipodal set in $(X,\| \cdot \|)$. It is easy to see that if
$\lambda>0$ then S is also bounded and separated with constants
$c_1,\lambda c_2,\lambda d$ and the same is valid for the set
$\lambda S=\{\lambda x:x \in S\}$ with constants $\lambda
c_1,c_2,\lambda d$.

(2) It follows from the above remark that an antipodal bounded and
separated set can be defined as a subset S of $B_X$ satisfying the
property that there is a constant $d>0$ such that for every $x,y
\in S$ with $x \neq y$ there exists $f \in B_{X^*}$ with $d \le
f(y)-f(x)$ and $f(x) \le f(z) \le f(y)$ for $z \in S$; that is, we
may assume that $c_1=c_2=1$. Given that formulation of Definition
2, it would be interesting to know if every \ind space contains an
infinite antipodal bounded and separated set with ($c_1=c_2=1$
and) $d>1$. (The answer is positive in case when X contains some
$\ell_p, 1 \le p <\infty$ or $c_0$). We note in this connection
that by a result of Elton and Odell the unit sphere of every \ind
space contains an infinite $(1+\varepsilon)$-separated set for
some
$\varepsilon>0$ (\cite{EO}, see also \cite{D} and \cite{KP}).\\

\noindent \textbf{Examples} Let X be a Banach space.

(1) Each finite antipodal set in X is bounded and separated
(obvious).

(2) Let $\{(x_\gamma,x_\gamma^*),\gamma \in \Gamma \}$ be a
bounded biorthogonal system in X; that is, there is a constant
$M>0$ such that $\|x_\gamma\| \cdot \|x_\gamma^*\| \le M$ for all
$\gamma \in \Gamma$. We set
$y_\gamma=\frac{x_\gamma}{\|x_\gamma\|}$ and
$y_\gamma^*=\|x_\gamma\| \cdot x_\gamma^*$, for $\gamma \in
\Gamma$. Clearly the system $\{(y_\gamma,y_\gamma^*),\gamma \in
\Gamma \}$ is biorthogonal. Now it is easy to see that the minimal
system $\{y_\gamma:\gamma \in \Gamma\}$ is (by Prop.1) antipodal
bounded and separated, with constants $c_1=1,c_2=M$ and $d=1$.

(3) Each \eq set S in X is bounded and separated antipodal set.
Indeed, as it follows from the method of proof of Prop.2, if S is
$\lambda$-\eq then the desired constants are $c_1=M, c_2=1$ and
$d=\lambda$, where $M=sup\{\|x\|:x \in S\}$. (Each \eq set is
clearly bounded).

The following result generalizes simultaneously a result of Petty
(\cite{Pe}, Th.2) and a result of Swanepoel already mentioned in
the introduction (\cite{S1}).

\begin{theo}
Let $(X,\| \cdot \|)$ be a \ban and $S \subseteq X$ be a bounded
and separated antipodal set. Then we have:
\begin{enumerate}
\item There is an equivalent norm $||| \cdot |||$ on $X$, such that
S is an \eq set in $(X,||| \cdot |||)$.
\item If the constants of S are $c_1=1,c_2=c$ and d, then the
Banach-Mazur distance between $(X,\| \cdot \|)$ and $(X,||| \cdot
|||)$ satisfies the inequality\\ $d((X,\| \cdot \|),(X,||| \cdot
|||)) \leq 2 \cdot \frac{c}{d}$.
\end{enumerate}
\end{theo}

\beg Assume (as we may) that the constants of S are 1,c and d (see
Remarks 3). We set
\[K=\overline{conv}(\frac{d}{c} \cdot B_X \cup \{x-y:x,y \in S \}).\] Then K
is a closed (bounded), convex symmetric set with $0 \in int(K)$,
so the corresponding Minkowski functional defines a norm on X
\[\|x\|_K=inf\{\lambda>0: x \in \lambda K\}\] and the unit ball of the space
$(X,\| \cdot \|_K)$ is exactly the set K. For $x,y \in S,x \neq y$
there is an $f \in cB_{X^*}$ such that:
\[d \leq f(y)-f(x) \leq \|f\| \|x-y\| \leq 2c\] hence $\frac{d}{c} \cdot B_X
\subseteq K \subseteq 2 \cdot B_X$ so it follows that the
Banach-Mazur distance of the two norms is $\leq 2 \cdot
\frac{c}{d}$.

It suffices to show that, if $x,y \in S$ with $x \neq y$, then
$x-y \in \partial K$ (equivalently $\|x-y\|_K=1$, where $\partial
K$ stands for the boundary of the set K) from which we have that S
is a 1-\eq set in $(X,\|\cdot \|_K)$. Let $x,y \in S$ with $x \neq
y$. Then there is $f \in cB_{X^*}$ with $d \leq f(y)-f(x)$ and
$f(x) \leq f(z) \leq f(y) \, \forall z \in S$. For every $z_1,z_2
\in S$ we have $f(z_1-z_2) \leq f(y-x)$. Also if $z \in
\frac{d}{c} \cdot B_X$, then $f(z) \leq |f(z)| \leq \|f\| \|z\|
\leq d \leq f(y-x)$, hence f is a support functional of the set K
through the point $y-x$ and so $y-x \in \partial K$.
\end{proof}

\begin{cor}
Let $\{(x_\gamma,{x_\gamma}^*):\gamma \in \Gamma \}$ be a bounded
biorthogonal system in the \ban $(X,\| \cdot \|)$, such that
$\|x_\gamma\|=1$ and $\|x_\gamma^*\| \le c$ for all $\gamma \in
\Gamma$. Then there is an equivalent norm $||| \cdot|||$ on X of
Banach-Mazur distance at most 2c from the original norm, such that
$\{x_\gamma:\gamma \in \Gamma\}$ is an \eq set in $(X,||| \cdot
|||)$.
\end{cor}

\beg The set $\{x_\gamma:\gamma \in \Gamma\}$ is bounded and
separated antipodal set with constants $c_1=1,c_2=c$ and $d=1$. So
theorem 3 can be applied.
\end{proof}

\noindent \textbf{Remarks 4} (1) For a Banach space
$(X,\|\cdot\|)$ with $dimX=\infty$ set $ant(X)=\sup\{d>0:$ there
is an infinite antipodal, bounded and separated set $S \subseteq
X$ with constants $c_1=c_2=1$ and $d\}$. We note that $ant(X) \le
2$; since by a result of Day (Th.1.20 in \cite{HMVZ}) there is an
infinite Auerbach system $\{(x_n,x_n^*): n \ge 1\}$ in X, that is,
a biorthogonal system with $\|x_n\|=\|x_n^*\|=1$ for $n \in
\Bbb{N}$ we get that $ant(X) \ge 1$. So Theorem 3 yields that for
every $\varepsilon>0$ X admits an equivalent norm with
Banach-Mazur distance $\le \frac{2}{ant(X)-\varepsilon}$ from the
original one, admitting an infinite equilateral set (c.f. Theorem
1 and Remark 3).

(2) A concept weaker than biorthogonality is that of
semibiorthogonality. Let X be a Banach space. A family
$\{(x_\alpha,x_\alpha^*):\alpha<\omega_1 \}$ is said to be
$\omega_1$-semibior\-thogonal, if it satisfies the following: (i)
$x_\beta^*(x_\alpha)=0$ for $\alpha<\beta<\omega_1$, (ii)
$x_\beta^*(x_\beta)=1$ for $\beta<\omega_1$ and (iii)
$x_\beta^*(x_\alpha) \ge 0$ for $\beta<\alpha<\omega_1$.

If we replace condition (iii) by the stronger: (iv) $0 \le
x_\beta^*(x_\alpha) \le 1$ for all $\alpha,\beta<\omega_1$, and if
the sets $\{x_\alpha:\alpha<\omega_1\}$
,$\{x_\alpha^*:\alpha<\omega_1\}$ are bounded then it is easy to
see that $\{x_\alpha:\alpha<\omega_1\}$ is a bounded, separated
and antipodal set. If for instance the compact space K contains a
closed non-$G_\delta$ set, then the \ban $C(K)$ admits an
$\omega_1$-semibior\-thogonal system of the form
$\{(f_\alpha,\delta_{t_\alpha}):\alpha<\omega_1\}$, where
$f_\alpha:K \rightarrow [0,1],\alpha<\omega_1$ are continuous
functions and $\delta_t$ is the Dirac measure at $t \in K$ (see
\cite{HMVZ} Prop. 8.7.). The compact scattered non-metrizable
space constructed (under CH) by Kunen is such that the \ban $C(K)$
admits an $\omega_1$-semibior\-thogonal system but no uncountable
biorthogonal system (\cite{HMVZ}, Th. 8.8 and Th. 4.41). We also
note that it is consistent with ZFC that there exist nonseparable
Banach spaces (of the form $C(K)$, where K is compact) which admit
no $\omega_1$-semibiorthogonal system (see \cite{LAT}, \cite{BGT}
and \cite{K}).

(3) It should be mentioned that there exist several interesting
classes of nonseparable Banach spaces, such as weakly compactly
generated (WCG) and their generalizations, that admit uncountable
bounded biorthogonal systems (actually Markushevich bases), see
\cite{HMVZ}. Finally notice that by a result of Todorcevic it is
consistent with ZFC (under Martin's Maximum axiom) to assume that
every nonseparable \ban admits an uncountable bounded biorthogonal
system (see \cite{To} and \cite{HMVZ}, Th. 4.48 and 8.12).

\noindent S.K.Mercourakis, G.Vassiliadis\\
University of Athens\\
Department of Mathematics\\
15784 Athens, Greece\\
e-mail: smercour@math.uoa.gr

\hspace{0.3cm} georgevassil@hotmail.com
\end{document}